\documentclass[12pt,a4paper,reqno]{amsart}%{article}

\def\a{{\bold{a}}}

\usepackage{amsmath}
\usepackage{amssymb}
\usepackage{amsfonts}

\advance\textheight by \topskip
\def\a{\alpha}
\def\b{\beta}

\title[On a functional--difference equation]
{On a functional--difference equation of Runyon, Morrison, 
Carlitz, and Riordan}

\author{Helmut Prodinger}

\address{ Helmut Prodinger,
The John Knopfmacher
Centre for Applicable Analysis and Number Theory,
School of Mathematics,
University of the Witwatersrand, P.\;O. Wits, 
2050 Johannesburg, South Africa, email:
{\tt helmut@gauss.cam.wits.ac.za},
 { www--address:{
}}{\tt http://www.wits.ac.za/helmut/index.htm}\,.
}

\keywords{functional-difference equation, generating function, 
Lagrange inverion formula, kernel method}
\subjclass{Primary:  05A15}

\date{\today}

\begin{document}
\begin{abstract}A certain functional--difference equation
that Runyon encountered when analyzing a queuing system
was solved in a combined effort of Morrison, Carlitz,
and Riordan.
We simplify that analysis by exclusively
using generating functions, in particular
 the {\em kernel method,\/} and the {\em Lagrange inversion
formula.}
\end{abstract}

\maketitle

\section{The equation}

The functional--difference equation in the title is
\begin{equation}\label{eq}   
(x-\a)(\a-\b)^{n-1}g_n(x)=\a(x-\b)^ng_{n-1}(\a)-x(\a-\b)^ng_{n-1}(x)
,\quad n\ge 1, \ g_0(x)=1.
\end{equation}

The aim of this note is to present a (possibly) simpler
solution than the (combined) solution by Morrison,
Carlitz, and Riordan \cite{Morrison64, Carlitz64, Riordan66}.

We introduce the generating function 
\begin{equation*}   
G(t):=\sum_{n\ge0}(\a-\b)^{n-1}g_n(x)t^n.
\end{equation*}
Multiplying (\ref{eq}) by $t^n$ and summing we get
\begin{equation}\label{erz}   
G(t)=\frac{\a\sum_{n\ge1}(x-\b)^nt^ng_{n-1}(\a)+
\frac{x-\a}{\a-\b}}{x-\a+xt(\a-\b)^2}.
\end{equation}
Now for 
\begin{equation*}   
x=\bar x=\frac{\a}{1+t(\a-\b)^2}
\end{equation*}
the denominator of (\ref{erz}) vanishes. Consequently,
the numerator must also vanish. (A more elaborate
argument would be that the power series expansion
must exist for that combination of values.)
This  is reminiscent of Knuth's trick
\cite[page 537]{Knuth97}, which is called {\em kernel method\/}
by some french authors. It leads to 
\begin{equation*}   
\sum_{n\ge1}(\bar x-\b)^nt^ng_{n-1}(\a)=
\frac{\bar x-\a}{(\b-\a)\a}.
\end{equation*}
Now we set $T=(\bar x-\b)t$, i.~e.
\begin{equation*}   
t=\frac{1-T(\a-\b)-\sqrt{1-2T(\a+\b)+T^2(\a-\b)^2}}
{2\b(\a-\b)}.
\end{equation*} 
So
\begin{equation*}   
\sum_{n\ge0}T^ng_n(\a)=\frac{1+T(\a-\b)-
\sqrt{1-2T(\a+\b)+T^2(\a-\b)^2}}{2T\a}.
\end{equation*}

The expansion of this generating function is
well known, from  the context of the Narayana
(Runyon!) numbers \cite{Riordan68} or elsewhere. In any instance, 
the coefficients could be easily detected by
the Lagrange inversion formula, with the result
\begin{equation*}   
g_n(\a)=\frac{1}{n}\sum_{k=0}^{n-1}\binom{n}{k}
\binom{n}{k+1}\b^{n-k}\a^k,\quad n\ge1,\ g_0(\a)=1.
\end{equation*}
In the next section, we will see a more impressive
occurrence of the Lagrange inversion formula.

\section{The general case}

In this section we move from the particular case
of $g_n(\a)$ to the general case of $g_n(x)$.

Now that the series in the numerator of (\ref{erz})
is established, the generating function $G(t)$ is
fully explicit:
\begin{equation}\label{erzz}   
G(t)=\frac{
{1+t(x-\b)(\a-\b)-
\sqrt{1-2t(x-\b)(\a+\b)+t^2(x-\b)^2(\a-\b)^2}}
+\frac{2(x-\a)}{\a-\b}}{2\big(x-\a+xt(\a-\b)^2\big)},
\end{equation}
 and one could work out some clumsy
expressions for the coefficients, e.~g. (for $x\neq\a$) 
\begin{equation*}   
g_n(x)=
\frac{(\a-\b)^{n}x^n}{(\a-x)^n}-
\a\sum_{k=1}^nx^{n-k}(\a-x)^{k-1-n}(\a-\b)^{n+1-2k}
(x-\b)^kg_{k-1}(\a).
\end{equation*}
This was obtained by Morrison without using the generating function.
Carlitz \cite{Carlitz64} set
\begin{equation}\label{defrel}
g_n(x)=\sum_{k=0}^{n-1}A_k^{(n)}(\a-\b)^{-k}(x-\b)^{k}
\end{equation}
and managed  to express the coefficients as follows:

\begin{equation*}
A_r^{(n)}=\b\phi_{r,n-1}-\a\sum_{s=1}^{r-1}g_{r-s}(\a)\, \phi_{s-1,n-r+s-1}
-\b\phi_{r-1,n-1},
\end{equation*}
with
\begin{equation*}   
\phi_{r,k}=\sum_{j=0}^{\min\{r,k\}}\binom{r}{j}\binom{k}{j}
\a^j\b^{k-j},\quad k\ge0, \qquad \phi_{r,k}=0, \quad k<0.
\end{equation*}
He asked whether the expressions

\begin{equation*}   
\mathcal{C}_{r,n}:=\sum_{s=1}^{r-1}g_{r-s}(\a)\, \phi_{s-1,n-r+s-1}
\end{equation*} 
can be simplified. 
Now Riordan \cite{Riordan66} proved that
\begin{equation*}
A_{k}^{(n)}=(n-k)\sum_{j=1}^k\frac1j\binom{n-1}{j-1}
\binom{k-1}{j-1}\a^j\b^{n-j},\qquad 1\le k<n,
\end{equation*}
and $A_{0}^{(n)}=\b^n$. (This was then
generalized by Carlitz
\cite{Carlitz66} who produced a $q$--version of that.) 
Riordan's answer translates as
\begin{equation*}   
\mathcal{C}_{r,n}=\sum_{j=1}^{\min\{r,n\}}
\binom{\min\{r,n\}}{j}\binom{\max\{r,n\}-1}{j-1}
\a^{j-1}\b^{n-j}.
\end{equation*}

Plugging Riordan's formula into the defining relation (\ref{defrel}) we get
(setting $w=\frac{\b-x}{\b-\a}$)
\begin{align*}
g_n(x)&=\b^n+\sum_{k=1}^{n-1}A_k^{(n)}w^{k}\\
&=\b^n+\sum_{k=1}^{n}(n-k)\sum_{j=1}^k\frac1j\binom{n-1}{j-1}
\binom{k-1}{j-1}\a^j\b^{n-j}w^{k}\\
&=\b^n+\sum_{j=1}^{n}\frac1j\binom{n-1}{j-1}
\a^j\b^{n-j}
\sum_{k=j}^{n}(n-k)
\binom{k-1}{j-1}w^{k}\\
&=\b^n+\sum_{j=1}^{n}\frac1j\binom{n-1}{j-1}
\a^j\b^{n-j}
[z^{n-1}]\frac{(zw)^j}{(1-z)^2(1-zw)^j}\\
&=\b^n+[z^{n-1}]\frac 1{(1-z)^2}\frac1n\sum_{j=1}^{n}\binom{n}{j}
\a^j\b^{n-j}
\frac{(zw)^j}{(1-zw)^j}\\
&=\frac1n
[z^{n-1}]\frac 1{(1-z)^2}\left(\frac{(\a-\b) zw+\b}{1-zw}\right)^n.
\end{align*}
This form was not observed before\footnote {Maple V.\;4 computed the inner sum in the
third line incorrectly, which cost me several hours!}.
It reminds one of the Lagrange inversion formula!
Consequently (see e.~g. \cite{Stanley99, Wilf94} for the Lagrange inversion formula)
\begin{equation}\label{lif}
g_n(x)=[t^n]\frac1{1-y},\quad\text{where}\quad y=t\Phi(y),
\quad\text{with}\quad \Phi(y)=\frac{(\a-\b) yw+\b}{1-yw}.
\end{equation}
This can be made explicit, since
\begin{equation*}
y=\frac{\scriptstyle{\a-\b+t(\b-\a)(\b-x)-\sqrt{(\b-\a)^2
 -2 (\b-\a ) (\a+\b ) (\b-x )t+(\b-\a )^{2} (\b-x)^{2}{t}^{2}
}}}{2(x-\b)}.
\end{equation*}
Although in general, (\ref{lif}) holds only for $n\ge1$, it is valid here for
$n=0$ as well. Therefore we get the generating function
\begin{equation*}
\sum_{t\ge0}t^ng_n(x)=\frac1{1-y}.
\end{equation*}
This matches with the previous expression for $G(t)$ in (\ref{erzz}). Reading all the steps of 
this section backwards, we obtain a proof of Riordan's result, purely by 
the use of generating functions, avoiding any recursions and any guesswork
(as in \cite{Riordan66}).

\bibliographystyle{plain}

%\bibliography{C:/texfiles2001/pro_bib}

\end{document}